\documentclass[12pt]{amsart}
\usepackage{amsmath,amsthm}
\numberwithin{equation}{section}
\newtheorem{thm}[equation]{Theorem}

\newtheorem{lm}[equation]{Lemma}
\newtheorem{prp}[equation]{Proposition}

\theoremstyle{definition}

\theoremstyle{remark}
\newtheorem{rem}[equation]{Remark}
\newcommand{\sprf}{\noindent{\it Proof.}}
\newcommand{\sqed}{\hfill\rule{1.3mm}{3mm}\medskip}

\newcounter{stareq}
\def\thestareq{\fnsymbol{stareq}}

\DeclareMathOperator{\BZ}{\mathbb{Z}} 
\DeclareMathOperator{\BR}{\mathbb{R}} 

\DeclareMathOperator{\inter}{\cap}  

\newcommand{\bd}{\begin{description}}
\newcommand{\ed}{\end{description}}
\newcommand{\vep}{\varepsilon}

\begin{document}

\date{\today}

\title{The Boltzmann--Sinai Ergodic Hypothesis \\ In Full Generality}

\author{N\'andor Sim\'anyi}

\address{The University of Alabama at Birmingham \\
Department of Mathematics \\
1300 University Blvd., Suite 452 \\
Birmingham, AL 35294 U.S.A.}

\email{simanyi@math.uab.edu}

\thanks{Research supported by the National Science Foundation, grant DMS-0800538}



\begin{abstract}
In the ergodic theory of semi-dispersing billiards the Local Ergodic
Theorem, proved by Chernov and Sinai in \cite{SCH87}, plays a central
role. So far, all existing proofs of this theorem had to use an
annoying global hypothesis, namely the almost sure hyperbolicity of
singular orbits. (This is the so called Chernov--Sinai Ansatz.)  Here
we introduce some new geometric ideas to overcome this difficulty and
liberate the proof from the tyranny of the Ansatz. The presented proof
is a substantial generalization of my previous joint result with
N. Chernov \cite{CS10} (which is a $2D$ result) to arbitrary
dimensions.

An important corollary of the presented ansatz-free proof of the Local
Ergodic Theorem is finally completing the proof of the
Boltzmann--Sinai Ergodic Hypothesis for hard ball systems in full
generality.
\end{abstract}

\maketitle

\parskip=0pt plus 2.5pt

\section{Introduction} \label{introduction}

Semi-dispersing billiards constitute an important class of finite-dimensional 
classical mechanical models, helping us understand deep connections between macroscopic 
and microscopic behaviours in statistical physics and hydrodynamics. Being such, studying
the ergodic properties of these systems has long been a substantial part of the theory of
dynamical systems and mathematical physics.

The definition of semi-dispersing billiards, in a nutshell, is as follows. We consider a
compact domain (the so called billiard table) $\mathbf{Q}\subset\mathbb{T}^d=\mathbb{R}^d/\mathbb{Z}^d$
or $\mathbf{Q}\subset\mathbb{R}^d$ whose boundary $\partial\mathbf{Q}$ is a finite union of $C^2$-smooth
hypersurfaces, so that the pair $(\mathbf{Q},\, \partial\mathbf{Q})$ has a cell-complex structure, the smooth
hypersurface components $(\partial\mathbf{Q})_i$ of $\partial\mathbf{Q}$ are concave from the outside of
$\mathbf{Q}$, i. e. they are ``bending away'' from $\mathbf{Q}$, and $\text{Int}(\mathbf{Q})$ is connected.
In addition to the above conditions we always assume that the $d$-dimensional spatial angle subtended by
$\mathbf{Q}$ at any of its boundary points $q\in \partial\mathbf{Q}$ is positive. In light of
\cite{BFK98}, this condition guarantees that only finitely many collisions can occur in finite time on
any trajectory, hence the full orbit of any non-singular phase point is uniquely defined. The only exceptional
category of semi-dispersing billiards on which we do not impose this angle-positivity condition is
the class of hard ball systems. For these systems the finiteness of the number of collisions
(in finite time) is guaranteed by another result of \cite{BFK98}.
The semi-dispersing billiard flow $\left\{\Phi^t\right\}$ describes the inertia motion of a point particle
inside $\mathbf{Q}$, undergoing specular reflections (the angle of reflection equals the angle of
incidence) when hitting the boundary $\partial\mathbf{Q}$. The phase space $\mathbf{M}$ of the arising flow
$\left(\mathbf{M},\, \left\{\Phi^t\right\},\, \mu\right)$ is essentially $\mathbf{Q}\times\mathbb{S}^{d-1}$,
modulo the natural identifications of the incoming and outgoing velocities $v^-$ and $v^+$ at the
boundary $\partial\mathbf{Q}$, and the invariant, finite Liouville measure $\mu$ is the product of the
Lebesgue measures of $\mathbf{Q}$ and of the unit velocity sphere $\mathbb{S}^{d-1}$. 

From the viewpoint of abstract ergodic theory, the two most fundamental problems concerning semi-dispersing
billiard systems are the full hyperbolicity (i. e. nonzero relevant Lyapunov exponents $\mu$-almost
everywhere) and the ergodicity of $\left(\mathbf{M},\, \left\{\Phi^t\right\},\, \mu\right)$. According to
the results of Chernov-Haskell \cite{CH96} and Ornstein-Weiss \cite{OW98}, the Bernoulli property of
semi-dispersing billiards follows, once the full hyperbolicity and the ergodicity are proven.

The origin of semi-dispersing billiards goes back (at least) to
Sinai's seminal papers \cite{S63} and \cite{S70}, in which he defined
semi-dispersing billiards to a large extent in the context of a
smaller class of models, namely the so called hard ball systems in
flat tori, shortly hard ball systems. In this context he gave a modern
formalism of Ludwig Boltzmann's ``ancient'' ergodic hypothesis: Every
hard ball system on a flat torus is fully hyperbolic and ergodic,
after we fix the values of the trivially conseved invariant
quantities, that is, the total momentum (set to zero), the center of
mass, and the total kinetic energy (set to $1$).  Nowadays, this
problem or hypothesis is widely known by the name of 
``Boltzmann--Sinai Ergodic Hypothesis''.

The hard ball systems with $n$ elastically colliding balls on 
$\mathbb{T}^\nu=\mathbb{R}^\nu/\mathbb{Z}^\nu$ reduce to
the motion of a billiard particle in a $\nu (n-1)$-dimensional torus
bouncing off $n(n-1)/2$ cylindrical obstacles (the billiard particle
hits a cylinder whenever two balls collide). Billiards with
cylindrical walls belong to a more general category of semi-dispersing
billiards, where a particle moves in a container with concave (but not
necessarily strictly concave) boundaries.

We remark that in the case $n=2$ the cylinders actually become
spheres, i.e.\ any system of $2$ hard balls reduces to a billiard
particle on a torus with a spherical obstacle. Such billiards belong
to a more special class of dispersing billiards, where a particle
moves in a container with strictly concave walls.

Dispersing billiards are always completely hyperbolic and ergodic
\cite{SCH87}, but for semi-dispersing billiards this may not be true.
For example, a billiard in a 3-torus with a single cylindrical wall
has zero Lyapunov exponents and is not ergodic; on the other hand, 2
transversal cylindrical walls within a 3-torus ensure hyperbolicity
and ergodicity \cite{KSSZ89}. For the systems of $n \geq 3$ hard
balls, one has to carefully explore the geometry of the cylindrical
walls in order to derive hyperbolicity and ergodicity.

There are two complications in the study of hard balls (or more
generally, semi-dispersing billiards). One is caused by the
\emph{singularities} of the dynamics -- these happen during
simultaneous multiple collisions of $\geq 3$ balls and during
grazing (tangential) collisions. In the phase space, singular points
make submanifolds of codimension one. The other complication is
caused by \emph{non-hyperbolicity} (i.e.\ the existence of zero
Lyapunov exponents) at some phase points. Such points make various
structures, ranging from  smooth submanifolds to Cantor-like subsets
of the phase space.

Powerful techniques have been developed to handle these two
complications separately (singularities and non-hyperbolicity), but
the combination of the two still presents a barely manageable situation.
More precisely, if non-hyperbolic sets and singularities intersect
in a subset of positive $[2\nu(n-1)-2]$-dimensional measure, then
modern proofs of ergodicity stall. On the other hand, such
substantial overlaps between singularities and non-hyperbolic sets
appear very unlikely (physically); they are regarded as
``conspiracy'', a word coined by M. Wojtkowski.

To bypass this scenario in an early work, Ya.~Sinai and N.~Chernov
\cite{SCH87} \emph{assumed} that almost every point on the
singularity manifolds (with respect to the intrinsic Lebesgue
measure) was completely hyperbolic (sufficient). Under this assumption (now
referred to as Chernov-Sinai Ansatz) they proved the so-called Local
Ergodic Theorem (also called `Fundamental Theorem'), which later
became instrumental in the proofs of ergodicity for various
billiards \cite{KSSZ90}, \cite{LW95}, \cite{BCST02}. It gives 
sufficient (and easily verifiable) conditions under which a phase 
point has an open neighborhood which belongs (mod 0) to one ergodic 
component.

Since \cite{S63} there have been several attempts to prove the full hyperbolicity
and/or the ergodicity of different hard ball systems. It is not our goal to present here
a comprehensive list of such results; what follows is only a (perhaps, subjective) selection
of some of them.

A.~Kr\'amli, N.~Sim\'anyi and D.~Sz\'asz built upon the results of
\cite{SCH87} and established the ergodicity for systems of $n=3$ hard
balls in any dimension \cite{KSSZ91}, and for $n=4$ hard balls in
dimension $d\geq 3$ \cite{KSSZ92}; in particular they verified
Chernov-Sinai Ansatz in these cases. However, their techniques could
not be extended to $n \geq 5$. The situation called for novel
approaches.

A partial breakthrough was made by Sim\'anyi and Sz\'asz when they
invoked ideas of algebraic geometry to rule out various
``conspiracies'' (at least for generic systems of hard balls), which
were in the way of proving hyperbolicity and ergodicity. Precisely,
they assumed that the balls had arbitrary masses $m_1,\dots,m_n$
(but the same radius $r$) and proved \cite{SSZ99} complete
hyperbolicity for generic vectors of ``external parameters'' 
$(m_1, \ldots, m_n, r)$; the latter needed to avoid some
exceptional submanifolds of codimension one in $\mathbb{R}^{n+1}$, which
remained unspecified and unknown. Later in \cite{S03} and \cite{S04} 
I used the same approach to prove Chernov-Sinai Ansatz
and ergodicity for generic systems of hard balls (in the above
sense). I also established hyperbolicity for systems of hard balls
of arbitrary masses \cite{S02}.

Thus the Boltzmann-Sinai ergodic hypothesis is now proved for
typical, or generic, systems of hard balls. This seems to be a
comforting settlement in both topological and measure-theoretic
senses, but it falls short of solving physically relevant problems,
as there is no way to check whether any particular system of hard
balls is ergodic or not. Most notably, for the system of balls with
all equal masses (which lies in the foundation of statistical
mechanics) the ergodicity remained open.

In an attempt to extend his results to \emph{all} gases of hard balls
(without exceptions), I developed \cite{S09} a new approach
based on purely dynamical (rather than algebro-geometric) ideas;
this allowed me to derive ergodicity from Chernov-Sinai Ansatz for
all hard ball systems. Thus the Boltzmann-Sinai hypothesis was
solved conditionally, modulo the Ansatz. It remained to prove the Ansatz,
or alternatively, derive the Local Ergodic Theorem without the Ansatz. (As a
side remark, it is ironic that the Ansatz, which originally seemed to be
just a convenient and temporary technical assumption, remained
the last unresolved issue in the whole picture.)

In \cite{CS10} Chernov and myself made another step toward a final
solution of the Boltzmann--Sinai Ergodic Hypothesis: we derived the Local
Ergodic Theorem without the Ansatz for arbitrary semi-dispersing
billiards in dimension $d=2$ and with algebraically defined smooth
components $(\partial\mathbf{Q})_i$ of the boundary $\partial\mathbf{Q}$.

\medskip

\begin{rem}
About the necessity of assuming the algebraic nature of $\partial\mathbf{Q}$ please
see \cite{BCST02}.
\end{rem}

\medskip

In this paper we prove the generalization of the result of \cite{CS10} to arbitrary
dimensions $d\geq 2$. For the definition of a sufficient (geometrically hyperbolic) orbit
$\Phi^{(-\infty,\infty)}x$ and for the Chernov--Sinai Ansatz, the reader is kindly referred to
Section 2 of \cite{S09}. In technical terms, our present result is as follows:

\medskip

\begin{thm} \label{main_result1}
Assume that the semi-dispersing billiard system $(\mathbf{M},\, \{\Phi^t\},\, \mu)$ 
has an algebraically defined boundary
$\partial\mathbf{Q}$.  Let $x_0 \in \partial\mathbf{M}$ be a point
whose entire trajectory $\Phi^{(-\infty, \infty)} (x)$ passes through
at most one singularity and is sufficient. Then there exists an open
neighborhood $U_0\subset\partial\mathbf{M}$ of $x_0$ in $\partial\mathbf{M}$ 
that belongs (mod 0) to one ergodic component of the billiard map $T$.
\end{thm}

\bigskip

We note that the base neighborhood $U_0=U_0(x_0)$ in this theorem is a
ball-like open neighborhood $U_0$ of $x_0$ consisting of sufficient
points for which

\medskip

(i) $U_0$ is a subset of the neighborhood $U_{\vep_1}(x_0)$ of $x_0$
featuring Theorem 3.6 in \cite{KSSZ90}, where $0<\vep_1<1$ is any fixed
number, and

\medskip

(ii) $U_0$ admits a family

$$
\mathcal{G}^\delta=\left\{G_i^\delta \big|\; i=1,2,\dots,I(\delta)
\right\} \quad (0<\delta<\delta_0)
$$
of regular coverings with a small enough threshold $\delta_0>0$, 
as explained in \cite{KSSZ90}.

\bigskip

All the existing proofs of the Local Ergodic Theorem \cite{SCH87},
\cite{KSSZ90} assume the Ansatz, and we relax that assumption.

Given a particular semi-dispersing billiard with an algebraically
defined $\partial\mathbf{Q}$, one can verify its ergodicity by showing
that the set of sufficient points is connected and has full measure.

\medskip

In light of the result of \cite{S09} and Theorem \ref{main_result1} we obtain the
following, important corollary.

\medskip

\begin{thm}[Verification of the Boltzmann--Sinai Ergodic Hypothesis] \label{main_result2}
Every hard ball system on a flat torus is fully hyperbolic and ergodic.
\end{thm}

\noindent
(For the definition of hard ball systems, please see, for example, Section 2 of 
\cite{S09}.)

\medskip

In proving Theorem \ref{main_result1}, in the subsequent sections we will be following the notations
of \cite{SCH87}, \cite{KSSZ90}, and \cite{CS10} as closely, as we can. 

The underlying idea of the present proof is borrowed from \cite{CS10} but, in order to adapt it to
the high-dimensional set-up, it had to undergo substantial improvements. This idea, in nutshell, is
as follows.

\medskip

For typical phase points $y\in U_0$, taken from a suitably chosen, small, open neighborhood
$U_0=U_0(x_0)$ of a sufficient phase point $x_0$, the local, stable invariant manifold
$\Sigma^s(y)$ of $y$ is to be constructed as the limiting manifold
$\lim_{t\to\infty} \Sigma^t_0(y)$, where $\Sigma^t_0(y)=\Phi^{-t}\left(\Sigma^t_t(y)\right)$,
$\Sigma^t_t(y)\ni\Phi^t y$ being a flat, maximal, local orthogonal manifold containing the 
phase point $y_t=\Phi^t y$. More precisely, $\Sigma^t_0(y)$ is the smooth component of the pre-image
$\Phi^{-t}\left(\Sigma^t_t(y)\right)$. The main problem in this construction arises when the boundary
$\partial\Sigma^t_0(y)$ of $\Sigma^t_0(y)$ gets ``very close'' to $y$, i. e. $\Sigma^t_0(y)$
is of ``very small size'' or, in other words, it has a very small internal radius around the center
$y$. This bad event is caused by the phenomenon that for some number $\tau$, $0\leq\tau\leq t$,
the phase point $y_\tau=\Phi^\tau y$ gets ``abnormally close'' to a singular reflection point
$w\in \mathcal{S}_0\subset\partial\mathbf{M}$. Of course, the concept of ``abnormally close''
is exactly defined in sections $3$--$4$. The above described phenomenon is the ``bad event'' that occurs
for the ``bad'' phase points $y\in U_0$. The proof of the Local Ergodic Theorem requires a strong, useful
upper measure estimate, the so called \emph{Tail Bound}, for the set of bad phase points $y\in U_0$.
In the original Chernov--Sinai approach this estimate was obtained by utilizing a hypothesis
(nowadays called Chernov--Sinai Ansatz), stating that almost every singular reflection point
$w'\in\mathcal{S}_0$ (with respect to the hypersurface measure of $\mathcal{S}_0$) is sufficient, i. e.
geometrically hyperbolic. If we eliminate this annoying hypothesis, we can still proceed with the
proof of the Local Ergodic Theorem by verifying that

\medskip

\begin{enumerate}
\item[(i)] a suitable approximating phase point $w'\in\mathcal{S}_0$ of $y_\tau$ is sufficient, and
\item[(ii)] $y_\tau$ and $w'$ can be joined by a flat hyperplane $\mathcal{P}$ that is
perpendicular to $\mathcal{S}_0$, unlike the curved hypersurface 
$\Sigma^t_\tau (y)=\Phi^{\tau-t}\left(\Sigma^t_t(y)\right)$ joining $y_\tau$ with
the original singular reflection point $w\in\mathcal{S}_0$. 
\end{enumerate}

\medskip

The original approximating singular reflection point 
$w\in\mathcal{S}_0\inter \Sigma^t_\tau (y)$ is automatically sufficient, for its pre-image
$\Phi^{-\tau} w$ happens to belong to the neighborhood $U_0$, this set containing only sufficient
phase points. When one replaces the pair $\left(\Sigma^t_\tau (y),\, w\right)$ with the pair
$(\mathcal{P},\, w')$, one needs a mechanism to guarantee the sufficiency of $w'$. The only imaginable
way to achieve this is the membership relation $\Phi^{-\tau} w'\in U_0$. The geometric ideas of Section 2
make sure that the ``tube'' 
\[
\bigcup_{0\leq s\leq \tau} \Phi^{-s}(\mathcal{P})
\]
is contained by the similarly formed ``tube'' 
\[
\bigcup_{0\leq s\leq \tau} \Phi^{-s}\left(\Sigma^t_\tau (y)\right).
\]
(This is Proposition \ref{main_lemma}.)

The subsequent section explains the role of the Tail Bound in the proof of the Local Ergodic Theorem,
along with pointing out that in the earlier proofs the Ansatz was only used to prove this bound.
Its main result is the ``inclusion lemma'', Lemma \ref{inclusion}, which appears to be a subtle
consequence of Proposition \ref{main_lemma}.

Finally, the closing section presents the entire, modified proof of the Tail Bound (along the lines
roughly sketched above) that spares the need to use the Chernov--Sinai Ansatz.

\bigskip

\section{Geometric Lemmas \\ Convex Hypersurfaces} \label{geometry}

In this section we will be dealing with {\it convex hypersurfaces} $\Sigma$ of the configuration
space $\mathbf{Q}$ of the considered semi-dispersing billiard flow. By a convex hypersurface
$\Sigma\subset\mathbf{Q}$ we shall mean a smooth (always meaning $C^\infty$ smooth) hypersurface
$\Sigma\subset\mathbf{Q}$ with a boundary $\partial\Sigma$ satisfying the following hypotheses:

\begin{enumerate}
\item[$(1)$] The pair $(\Sigma, \partial\Sigma)$ is diffeomorphic to the standard pair 
$(B^{d-1}, \mathbb{S}^{d-2})$, where $B^{d-1}$ is the closed unit ball of $\BR^{d-1}$ with the
boundary $\mathbb{S}^{d-2}=\partial B^{d-1}$.
\item[$(2)$] $\Sigma$ is equipped with a unit normal vector field $n(q)$ ($q\in\Sigma$), and the
second fundamental form $B_\Sigma(q)$ of $\Sigma$ at any point $q\in\Sigma$ is positive semi-definite.
\item[$(3)$] The boundary $\partial \Sigma$ of $\Sigma$ is geodesically strictly convex in 
$\Sigma$ in the sense that the second fundamental form of $\partial \Sigma$ inside 
$\Sigma$ (with respect to an outer, init normal vector field) is strictly positive at every
point of $\partial \Sigma$.
\end{enumerate}

The convex hypersurface $\Sigma\subset\mathbf{Q}$ supplied with the specified normal vector field
$n(q)$ ($q\in\Sigma$) shall always be denoted by $\tilde{\Sigma}$, i. e. 
\[
\tilde{\Sigma}=\left\{(q, n(q))\big|\; q\in\Sigma\right\}\subset\mathbf{M}.
\]

\medskip

\begin{rem}
According to the above definition, a convex hypersurface is convex in two different
senses: One of them is postulated in (2), whereas the other one in (3). The meaning 
of the latter one is that the boundary $\partial\Sigma$ of $\Sigma$ is bending toward $\Sigma$.
\end{rem}

\begin{rem}
Sometimes the manifold $\tilde{\Sigma} \subset \mathbf{M}$ will be called a convex hypersurface
and its projection $\Sigma=\pi\left(\tilde{\Sigma}\right)\subset\mathbf{Q}$ the
{\it carrier} of $\tilde{\Sigma}$. However, this little ambiguity in the notations should never
cause any confusion.

In this paper $\pi:\; \mathbf{M}\to \mathbf{Q}$ always denotes the natural projection.
\end{rem}

\begin{rem}
The essential part of the (stable version of the) Local Ergodic
Theorem \cite{SCH87}, see also \cite{KSSZ90}, is the construction of a
sufficiently large amount of ``not too short'' local stable invariant
manifolds. Actually, it turns out that all of the so called {\it local
  orthogonal manifolds} (LOMs) featuring in those constructions appear
to be convex hypersurfaces $\Sigma$ of $\mathbf{Q}$, as defined
above. More precisely, in Lemma 4 of \cite{SCH87} (the fundamental
construction of local stable invariant manifolds $\Sigma^s(y)$ for
typical phase points $y\in U^g(\delta)$, see also the counterpart
result, Lemma 5.4 of \cite{KSSZ90}), instead of taking the convex,
orthogonal manifolds $\Sigma^n_0(y)=T^n\left[\Sigma^n_n(y)\right]$ (as
the $n$-step approximation of $-\Sigma^s(y)$), one can take the
largest open, geodesic round ball $\mathcal{B}^n_n(y)$ inside
$\Sigma^n_n(y)$, centered at the configuration component $\pi(T^n y)$
of $T^n y$, on which the iterated map $T^n$ is smooth, and operate
with the manifold
$\mathcal{B}^n_0(y)=T^n\left[\mathcal{B}^n_n(y)\right]$ instead of
$\Sigma^n_0(y)$. We note here that the relevant metric on the manifold
$\Sigma^n_n(y)$ is the orthogonal configuration displacement $\Vert
\delta q\Vert$, please see the beginning of Section \ref{tail_bound}.
The largest open, geodesic ball in $\Sigma^n_n(y)$ is understood
with respect to this metric.  Then the manifold $\mathcal{B}^n_n(y)$
is a convex hypersurface in $\mathbf{Q}$, according to our
definition above. Finally, the limiting set
\[
\lim_{n\to \infty} -\mathcal{B}^n_0(y)
\]
will be a suitably big, open submanifold of $\Sigma^s(y)$ to be constructed in the cited lemmas.
\end{rem}

\begin{rem}
A convex hypersurface $\Sigma\subset\mathbf{Q}$ is said to be flat if $B_\Sigma(q)\equiv 0$ on
$\Sigma$ or, equivalently, if $\Sigma$ lies in a hyperplane of $\mathbf{Q}$. Furthermore, $\Sigma$
shall be called strictly convex if $B_\Sigma(q)>0$ for all $q\in\Sigma$.
\end{rem}

For $t\in\mathbb{R}$ we introduce
\begin{equation}
\Sigma^t=\pi\left[\Phi^t(\tilde{\Sigma})\right],
\end{equation}
and for any interval $I\subset\mathbb{R}$
\begin{equation}
\Phi^I(\Sigma)=\bigcup_{t\in I} \Sigma^t.
\end{equation}

Our first result in this section is, essentially, a simple observation.

\begin{lm} \label{dilation}
Assume that $y_0=(q_0, n_0)\in\mathbf{M}\setminus\partial\mathbf{M}$, 
$0<\rho_0<z_{tub}(y_0)$,
\[
\begin{aligned}
\Sigma_0 &=\left\{q\in\mathbf{Q}\big|\; (q-q_0)\perp n_0,\; ||q-q_0||\leq \rho_0\right\}, \\
\tilde{\Sigma}_0 &=\left\{(q, n_0)\big|\; q\in\Sigma_0\right\},
\end{aligned}
\]
$T<0$,
\[
\text{Int}\left[\Phi^{[T,0]}(\Sigma_0)\right]\inter \partial\mathbf{Q}=\emptyset,
\]
$\Sigma_1\subset \Phi^{[T,0]}(\Sigma_0)$ is a strictly convex hypersurface with 
$y_0\in\tilde{\Sigma_1}$, and that for all $t\in [T,0]$ the hypersurface
$\pi\left(\Phi^t(\tilde{\Sigma_1})\right)=\Sigma^t_1$
does not get focused (in particular, it is still a
convex hypersurface) and it is a subset of $\text{Int}\mathbf{Q}$.

We claim that for every $y_1\in\tilde{\Sigma}_1$ and for every $t$ with 
$\max\{T,\, -\tau(-y_1)\}<t<0$ the inequality
\[
\text{dist}\left(q+tn_0,\, \pi(\Phi^t y_1)\right)<\rho_0
\]
holds true. In particular, $\Phi^t(\Sigma_1)\inter \partial\mathbf{Q}=\emptyset$ for all $t$ with
\[
\max\left\{T,\, \sup\left\{-\tau(-y_1)\big|\; y_1\in
\tilde{\Sigma}_1\right\}\right\}<t<0.
\]
\end{lm}

\begin{sprf}
The statement directly follows from the obvious fact that the distance function 
$g(t)=\left\Vert q_0+tn_0-\pi(\Phi^t y_1)\right\Vert$ is strictly increasing in 
$t$ for all $t$ with 
\[
\max\{T,\, -\tau(-y_1)\}<t\leq 0, 
\]
due to the strict convexity of each $\Phi^t\left(\Sigma_1\right)$.
\end{sprf} \qed

\medskip

The next result claims that the forward tube $\Phi^{[0,T]}(\Sigma_0)$ of a flat, convex hypersurface
$\Sigma_0$ cannot escape from the similarly formed tube $\Phi^{[0,T]}(\Sigma_1)$ of a strictly convex
hypersurface $\Sigma_1$, once $\Sigma_0\subset \Phi^{[0,T]}(\Sigma_1)$.

\medskip

\begin{prp} \label{main_lemma}
Assume that $y_0=(q_0,n_0)\in \mathbf{M}\setminus \partial\mathbf{M}$, $\Sigma_0$ and $\Sigma_1$
are convex hypersurfaces in $\mathbf{Q}$, $\Sigma_0$ is flat, $\Sigma_1$ is strictly convex,
$y_0\in\text{Int}(\tilde{\Sigma}_0)\inter\text{Int}(\tilde{\Sigma}_1)$, $T>0$, $\Phi^T$ is smooth
on $\tilde{\Sigma}_1$, and $\Sigma_0 \subset \Phi^{[0,T]}(\Sigma_1)$. We claim that
\[
\Phi^{[0,T_1]}(\Sigma_0)\subset\Phi^{[0,T]}(\Sigma_1)
\]
for every number $T_1$ with $0<T_1\leq T$ and
\begin{equation} \label{gap}
T-T_1\geq \sup_{q\in\Sigma_0} d(q,\, \Sigma_1):=\delta_0.
\end{equation}
In particular, $\Phi^{T_1}$ is smooth on $\tilde{\Sigma}_0$.
\end{prp}

\begin{sprf}
Let $H \subset \BR^d$ be the hyperplane containing $\Sigma_0$. By extending $\Sigma_0$ (if necessary), we may
assume that $\Sigma_0 = H \cap \Phi^{[0,T]}(\Sigma_1)$ with the boundary 
$\partial\Sigma_0 = H \cap \partial\Phi^{[0,T]}(\Sigma_1)$. By continuity it is enough to show for every
\textit{interior} point $y_0'=(q_0',n_0)\in \text{Int}(\tilde{\Sigma}_0)$ that $\Phi^{[0,T_1]}y_0'$ stays inside
the ``tube'' $\mathcal{T} = \Phi^{[0,T]}(\Sigma_1)$. Thanks to the condition $T-T_1 \geq \delta_0$, the orbit segment
$\Phi^{[0,T_1]}y_0'$ could only escape from the tube $\mathcal{T}$ through its ``side'' 
$\Phi^{[0,T]}(\partial\Sigma_1)$. In order to prove the proposition we shall show that
\[
\Phi^{[0,T_1]}y_0' \cap \Phi^{[0,T]}(\partial\Sigma_1) = \emptyset.
\]

Assume, on the contrary, that there exists a time $t_0$, $0 < t_0 \leq T_1$, such that
$\Phi^{t_0} y_0'=y_1=(q_1,v_1)$ and $q_1\in \Phi^{[0,T]}(\partial\tilde{\Sigma}_1)$. 
Let $y_2=(q_2,v_2)\in \partial\tilde{\Sigma}_1$, $t_1\in [0,T]$ such that 
\[
q\left(\Phi^{t_1}y_2\right) = q\left(\Phi^{t_0}y_0'\right) = q_1.
\]
(These objects are clearly unique.) We define the strictly convex hypersurface
$\Sigma_3 \subset \mathbf{Q}$ as follows:
\[
\tilde{\Sigma}_3 = SC_{-y_0'}\left[\Phi^{t_0}\left(\{q_1\}\times \mathbb{S}^{d-1}\right)
\right] \cap \mathcal{T}.
\]
Here $SC_x [A]$ stands for the smooth component of the set $A\subset \mathbf{M}$ containing
the point $x$ and $-y_0' = (q_0',\, -n_0)$. It is clear that
\[
\partial\tilde{\Sigma}_3 = SC_{-y_0'}\left[\Phi^{t_0}\left(\{q_1\}\times \mathbb{S}^{d-1}\right)
\right] \cap \Phi^{[0,T]}\left(\partial\Sigma_1\right).
\]
The strictly convex hypersurfaces $\Sigma_1$ and $\Sigma_3$ touch the hyperplane $H$ at the points
$q_0$ and $q_0'$, respectively, but they lie on opposite sides of $H$.

Let $[q_3,\, q_4]$ ($q_3\in \Sigma_3$, $q_4\in \Sigma_1$) be a shortest line segment connecting a point of
$\Sigma_3$ with a point of $\Sigma_1$. It follows from our convexity conditions that
$q_3\in \text{Int}(\Sigma_3)$. For any boundary point $z\in \partial\Sigma_1$ the entire set
$\text{Int}(\Sigma_3)$ lies in the interior of the side of the hyperplane 
$L = \mathcal{T}_z\Phi^{[-\epsilon, \epsilon]}\left(\partial\tilde{\Sigma}_1\right)$ 
that contains $\Sigma_1$, hence the other
end $q_4\in \Sigma_1$ of the length-minimizing segment $[q_3,\, q_4]$ cannot be a boundary point 
(of $\Sigma_1$) either. The relations $q_3\in \text{Int}(\Sigma_3)$, $q_4\in \text{Int}(\Sigma_1)$
imply that the segment $[q_3,\, q_4]$ is perpendicular to both $\Sigma_3$ and $\Sigma_1$. This implies,
in turn, that 
\begin{enumerate}
\item [$(i)$] $q\left(\Phi^{-t_0}(q_3,\, n_3)\right) = q_1$, where $n_3 = \frac{q_4-q_3}{\Vert q_4-q_3\Vert}$,
so that $(q_3,\, n_3)\in \tilde{\Sigma}_3$;
\item [$(ii)$] Inside the tube $\mathcal{T}$ the point $q_4$ is the unique point of $\Sigma_1$ lying closest
to the point $q_1$.
\end{enumerate}
(For (i) we note that the relation $q_3\not=q_4$ can be obviously assumed.)

\medskip

However, the fact $q\left(\Phi^{t_1}(q_2,\, v_2)\right) = q_1$ implies that, inside $\mathcal{T}$, the unique
point of $\Sigma_1$ lying closest to $q_1$ is the boundary point $q_2$, not the interior point $q_4$.
The obtained contradiction completes the proof of the proposition.
\end{sprf} \qed

\medskip

This result shall be applied in the sequel within the following circumstances:

\medskip

There will be given a sufficient phase point
$x_0\in\mathbf{M}\setminus\partial\mathbf{M}$ with a non-singular
forward orbit and a small, nicely shaped, open neighborhood $U_0$ (in
which a sufficiently large collection of ``not too short'' local
stable invariant manifolds are to be constructed, just as requested by
the proof of the Local Ergodic Theorem \cite{SCH87}, see also
\cite{KSSZ90}), so that, in addition to all the conditions formulated in 
Proposition \ref{main_lemma},

\begin{enumerate}
\item[$(a)$] Every $x\in U_0$ is sufficient;
\item[$(b)$] $\Phi^T\left(\tilde{\Sigma}_1\right)\subset U_0$.
\end{enumerate}

As a consequence of Proposition \ref{main_lemma}, we will get that 
$\Phi^{T_1}\left(\tilde{\Sigma}_0\right)\subset U_0$, hence for every point
$y_0'\in\tilde{\Sigma}_0$ the sufficiency of $y_0'$ will be guaranteed by the fact
$\Phi^{T_1} y_0'\in U_0$.

\bigskip

\section{Local Ergodic Theorem and \\ Tail Bound} \label{tail_bound}

Now we turn to the proof of the Local Ergodic Theorem without using the Ansatz. We closely follow the lines
and notations of \cite{KSSZ90} that presents a clear and complete proof of the Local Ergodic Theorem.

For the given sufficient phase point $x_0\in\text{Int}\mathbf{M}$ 
(with a non-singular forward orbit $\Phi^{(0,\infty)}x_0$) we consider a sufficiently small, open neighborhood 
$U_0=U_0(x_0)$ of $x_0$, just as described in Theorem 3.6 of \cite{KSSZ90}.

Given a convex hypersurface $\tilde{\Sigma}\subset\mathbf{M}$ with a carrier
$\Sigma=\pi(\tilde{\Sigma})$, we use the metric on it defined by the distance along the
orthogonal manifold $\Sigma$, and denote by $\Vert\, .\, \Vert=\Vert \delta q\Vert$ the 
corresponding norm on the tangent spaces $\mathcal{T}_y\tilde{\Sigma}$, $y\in\tilde{\Sigma}$.
For a convex hypersurface $\tilde{\Sigma}\subset\partial\mathbf{M}$ we use the norm and metric on the
corresponding flow-sync convex hypersurface of $\mathbf{M}$, constructed right after the considered
collision $y\in \tilde{\Sigma}\subset\partial\mathbf{M}$.

For any $x\in \tilde{\Sigma}\subset\partial\mathbf{M}$ denote by $D^n_{x,\tilde{\Sigma}}$ the derivative
(linearization) of the map $T^n$ (restricted to $\tilde{\Sigma}$) at $x$. Its norm
$\Vert D^n_{x,\tilde{\Sigma}}\Vert$ will always be taken with respect to the norm $\Vert\delta q\Vert$
described above. If $\tilde{\Sigma}\subset\partial\mathbf{M}$ is a convex hypersurface, then 
$\Vert D^n_{x,\tilde{\Sigma}}\Vert\geq 1$ for every $n\geq 1$. Denote
\begin{equation}
\kappa_{n,0}(x)=\inf_{\tilde{\Sigma}} \left\Vert\left(D^n_{-T^nx,\tilde{\Sigma}}\right)^{-1}\right\Vert^{-1},
\end{equation}
where the infimum is taken for all convex hypersurfaces $\tilde{\Sigma}$ through $-T^nx$. This quantity is
the minimum local expansion of convex hypersurfaces on their way from $-T^nx$ back to $-x$. (We note that the
infimum is actually attained at the flat hyperplane, cf. \cite{CM06}.)

Given $\delta>0$ we denote 
\begin{equation}
\kappa_{n,\delta}(x)=\inf_{\tilde{\Sigma}}\inf_{y\in\tilde{\Sigma}}
\left\Vert\left(D^n_{y,\tilde{\Sigma}}\right)^{-1}\right\Vert^{-1},
\end{equation}
where the infimum is taken for all convex hypersurfaces $\tilde{\Sigma}$ through $-T^nx$
for which $T^n$ is smooth on $\tilde{\Sigma}$ and
$\text{dist}\left(-x,\, \partial(T^n\tilde{\Sigma})\right)\leq\delta$. We observe that 
$1\leq \kappa_{n,\delta}(x)\leq \kappa_{n,0}(x)$, and both $\kappa_{n,\delta}(x)$
and $\kappa_{n,0}(x)$ are non-decreasing in $n$. 

For $x\in\partial\mathbf{M}$ we denote
\begin{equation}
z_{\text{tub}}(x)=\sup_{\tilde{\Sigma}} \left\{\text{dist}(x,\, \partial\tilde{\Sigma})\big|\;
T\text{ is smooth on } \tilde{\Sigma}\right\},
\end{equation}
where the supremum is taken for all flat hyperplanes $\tilde{\Sigma}\subset\partial\mathbf{M}$
through $x$. This number is the radius of the maximal tubular neighborhood of the billiard link
joining $x$ with $Tx$. Note that $z_{\text{tub}}(-Tx)=z_{\text{tub}}(x)$.

\medskip

We denote by $\tilde{\Sigma}^u(x)$ and $\tilde{\Sigma}^s(x)$ the (local) unstable and stable manifolds
through $x$; the former is a convex hypersurface and a latter one is a concave one. We also put
$r^\alpha(x)=\text{dist}(x,\, \partial\tilde{\Sigma}^\alpha(x))$ for $\alpha=u,\, s$. It is known, cf.
Lemma 5.4 of \cite{KSSZ90}, that for every semi-dispersing billiard system there exists a constant
$c_3>0$ such that if 
\[
x\in U^g=U^g(\delta)=\left\{y\in U_0\big|\; \forall \, n>0\; \;
z_{\text{tub}}(-T^ny)\geq\left(\kappa_{n,c_3\delta}(y)\right)^{-1}\cdot c_3\delta\right\},
\]
then $r^s(x)\geq c_3\delta$. Thus the points of $U^g=U^g(\delta)$ (``good set'') have stable manifolds of
order $\delta$. A similar property holds for unstable manifolds. The set of points with shorter stable
manifolds (the ``bad set'') must be carefully analyzed. We take
\begin{equation}
\begin{aligned}
U^b_n&=U^b_n(\delta)=\left\{y\in U_0 \big|\;
z_{\text{tub}}(-T^ny)<\left(\kappa_{n,c_3\delta}(y)\right)^{-1}\cdot c_3\delta\right\}, \\
U^b&=U^b(\delta)=U_0\setminus U^g=\bigcup_{n\geq 1} U^b_n.
\end{aligned}
\end{equation}

In the construction of the local stable invariant manifolds (i. e. from the viewpoint of the proof of the
Local Ergodic Theorem) the really relevant bad sets are the sets
\begin{equation} \label{bad_set}
\begin{aligned}
\hat{U}^b_n(\delta)&=\hat{U}^b_n=\Big\{x\in U_0\big|\; \Phi^{(0,\infty)}x \text{ is non-singular, }
\exists \text{ a strictly convex hypersurface } \tilde{\Sigma}, \\
& -T^nx\in\tilde{\Sigma},\, \partial\tilde{\Sigma}\cap\mathcal{S}_0=\emptyset, \, 
T^n \text{ and } T^{-1} \text{ are smooth on } \text{Int}(\tilde{\Sigma}), \\
& T^{-1}(\tilde{\Sigma}) \text{ is convex, } T^{-1} \text{ is not smooth at some point }
x'\in\partial\tilde{\Sigma}, \\
& \text{dist}(-T^nx,\, x')<\left(\kappa_{n,c_3\delta}(x)\right)^{-1}\cdot c_3\delta,\,
T^n(\tilde{\Sigma})\subset -U_0 \Big\}, \\
\hat{U}^b(\delta)&=\hat{U}^b=\bigcup_{n\geq 1} \hat{U}^b_n(\delta).
\end{aligned}
\end{equation}
For any $x\in\hat{U}^b_n(\delta)$ and $\tilde{\Sigma}$ as in (\ref{bad_set}) above, we denote
\begin{equation}
\begin{aligned}
z(T^n x,\, \tilde{\Sigma})&=\min_{x'} \big\{\text{dist}(-T^nx,\, x')\big|\; x'\in\partial\tilde{\Sigma}, \\
&T^{-1} \text{ is not smooth at } x'\big\}.
\end{aligned}
\end{equation}
Observe that this quantity is at most $\left(\kappa_{n,c_3\delta}(x)\right)^{-1}\cdot c_3\delta$.
We also note that 
\begin{equation}
z_{\text{tub}}(T^nx)\leq z(T^n x,\, \tilde{\Sigma}).
\end{equation}

A crucial fact in the proof of the Local Ergodic Theorem is the following

\medskip

\begin{thm}[Tail Bound]
For any function $F(\delta)\nearrow \infty$, as $\delta\searrow 0$, the set 
\[
\hat{U}^b_\omega=\hat{U}^b_\omega(\delta)=\bigcup_{n>F(\delta)} \hat{U}^b_n(\delta)
\]
has measure $\mu_1\left(\hat{U}^b_\omega(\delta)\right)=o(\delta)$, i. e. 
$\lim_{\delta\to 0}\delta^{-1}\cdot\mu_1\left(\hat{U}^b_\omega(\delta)\right)=0$.
\end{thm}

\medskip

In fact, the derivation of the Local Ergodic Theorem from the Tail Bound does not use the Ansatz
(cf. \S5 of \cite{KSSZ90}), so we will not repeat that part of the proof here. In what follows, we prove
the Tail Bound {\it without using the Ansatz}. We shall prove the upper measure estimate of the Tail Bound
for the larger sets
\begin{equation} \label{bad_tilde}
\begin{aligned}
\tilde{U}^b_n&=\tilde{U}^b_n(\delta)=\Big\{x\in U_0\big|\; \exists \vep_1,\,
0<\vep_1<\tau(T^nx), \\
& \exists\text{ a post-sufficient point } y\in\text{Int}\mathbf{M},\, Ty\in\mathcal{S}_0, \\
& v(y)=v(T^nx)=v\left(\Phi^{\vep_1}(T^nx)\right),\, \Delta q=\left[q(y)-
q\left(\Phi^{\vep_1}(T^nx)\right)\right]\perp v(T^nx), \\
& (\Delta q,\, 0)\perp \mathcal{T}_yJ,\, \Vert\Delta q\Vert<
\left(\kappa_{n,c_3\delta}(x)\right)^{-1}\cdot c_3\delta \Big\}, \\
\tilde{U}^b_{n,m}&=\tilde{U}^b_{n,m}(\delta)=\left\{x\in\tilde{U}^b_n\big|\;
\Lambda^m\leq \kappa_{n,c_3\delta}(x)<\Lambda^{m+1}\right\},
\end{aligned}
\end{equation}
see also the definition of the corresponding set $U^b_{n,m}$ right before Lemma 6.3 in
\cite{KSSZ90}. In the above definition (\ref{bad_tilde}) the set $J\subset\mathbf{M}$ is the
codimension-one, smooth submanifold (near $y$)
\[
J=\left\{y'\in\mathbf{M}\big|\; Ty'\in\mathcal{S}_0\right\}
\]
containing the point $y$. Also, in our terminology the phrase ``past-sufficient phase point
$y\in\text{Int}\mathbf{M}$'' means that the phase point $T^{-1}y\in\partial\mathbf{M}$ is sufficient,
i. e. $\Phi^{(-\infty,-\tau(-y))}y$ is sufficient.

\medskip

\begin{lm} \label{inclusion}
The inclusion
\[
\hat{U}^b_n(\delta)\subset\tilde{U}^b_n(\delta)
\]
holds true.
\end{lm}

\begin{sprf}
Let $x\in\hat{U}^b_n(\delta)$, a strictly convex hypersurface $\tilde{\Sigma}\ni -T^nx$,
and a point $x'\in\partial\tilde{\Sigma}$ be given, as stipulated in the definition of the set
$\hat{U}^b_n(\delta)$. We are going to use Proposition \ref{main_lemma} as follows: 
First of all, we take a time-sync version $\tilde{\Sigma}_1\subset\mathbf{M}$ of 
$\tilde{\Sigma}$ so that

\medskip

\begin{enumerate}
\item[$(1)$] $\tilde{\Sigma}_1\subset\mathbf{M}$ is a strictly convex hypersurface;
\item[$(2)$] The restriction of the natural projection $T:\; \mathbf{M}\to\partial\mathbf{M}$
to $\tilde{\Sigma}_1$ maps $\tilde{\Sigma}_1$ onto its image $T\left(\tilde{\Sigma}_1\right)=\tilde{\Sigma}$
in a diffeomorphic fashion.
\end{enumerate}
Denote this diffeomorphism by $\mathcal{D}:\; \tilde{\Sigma}_1 \xrightarrow{\approx} \tilde{\Sigma}$.

Let $x''=\mathcal{D}^{-1}(x')$ and
\[
\Phi^{-\vep_1}(-T^n x)=-\Phi^{\vep_1}(T^n x)=\mathcal{D}^{-1}(-T^n x).
\]
We may assume that $\tilde{\Sigma}$ (hence $\tilde{\Sigma}_1$, as well) is selected in such a way that 
\[
z(T^n x,\, \tilde{\Sigma})=\text{dist}(-T^n x,\, x')=
\text{dist}\left(-\Phi^{\vep_1}(T^n x),\, x''\right)
\]
(i. e. the minimum value in the definition of $z(T^n x,\, \tilde{\Sigma})$ is attained at $x'$), and that the
entire manifold (with boundary) $\Sigma_1\subset\mathbf{Q}$ is a closed, geodesic ball centered at
$q\left(-\Phi^{\vep_1}(T^n x)\right)$ having the radius $z(T^n x,\, \tilde{\Sigma})$. 

Finally, we take $y_0=-\Phi^{\vep_1}(T^n x)$, and
\begin{equation}
\begin{aligned}
\tilde{\Sigma}_0 & =\big\{y \in\mathbf{M}\big|\; v(y)=v(y_0),\; \left(q(y)-q(y_0)\right)\perp v(y_0), \\
& \Vert q(y)-q(y_0)\Vert\leq z_{\text{tub}}(T^n x)=z_{\text{tub}}(y_0)\big\}.
\end{aligned}
\end{equation}
We remind the reader that $z_{\text{tub}}(T^n x)\leq z(T^n x,\, \tilde{\Sigma})$. Let
$-y\in\partial\tilde{\Sigma}_0$ be a point for which $T^{-1}(-y)=-T(y)\in\mathcal{S}_0$.
According to the inequality $z_{\text{tub}}(T^n x)\leq z(T^n x,\, \tilde{\Sigma})$ and Lemma \ref{dilation}, 
all the hypotheses of Proposition \ref{main_lemma} are satisfied. Therefore, according to
Proposition \ref{main_lemma}, we have $T^{n+1}(-y)\in U_0$, thus the phase point $y$ is past-sufficient.
Finally, the chain of inequalities
\[
\text{dist}(y_0,\, y)\leq z(T^n x,\, \tilde{\Sigma})\leq
\left(\kappa_{n,c_3\delta}(x)\right)^{-1}\cdot c_3\delta
\]
shows that the phase point $y$ satisfies all conditions required in the definition of the set
$\tilde{U}^b_n(\delta)$, thus indeed $x\in \tilde{U}^b_n(\delta)$, finishing the proof of the
lemma.
\end{sprf} \qed

\medskip

Now, with the above constructions and lemma, we are ready to complete the proof of the 
Tail Bound (without the Ansatz) in the next section. It goes along the same lines as
in Section 6 of \cite{KSSZ90}, but at many points the argument needs (small) modifications
in order to avoid the need to use the Ansatz.

\bigskip

\section{Proof of the Tail Bound \\ (Without Ansatz)} \label{proof_tail_bound}
Having in mind Lemma \ref{inclusion}, we are going to prove a strengthened version of the Tail Bound.

\medskip

\begin{thm}[Tail Bound (Stronger version)] \label{strong_tail_bound}
For any function $F(\delta)$, such that $F(\delta)\nearrow \infty$ as $\delta\searrow \infty$,
we have
\[
\mu_1\left(\tilde{U}^b_\omega(\delta)\right)=o(\delta),
\]
i. e. $\lim_{\delta\to 0} \delta^{-1}\cdot \mu_1\left(\tilde{U}^b_\omega(\delta)\right)=0$,
where
\[
\tilde{U}^b_\omega=\tilde{U}^b_\omega(\delta)=\bigcup_{n>F(\delta)} \tilde{U}^b_n(\delta).
\]
\end{thm}

\medskip

\begin{sprf}
The proof uses, among other things, a simple arithmetic lemma, Lemma 6.2 of \cite{KSSZ90}.

\medskip

\begin{lm}[Lemma 6.2 of \cite{KSSZ90}] \label{arithmetic}
Assume that for any $\delta>0$, $a^\delta_{n,m}\geq 0$, $n,\, m\in\BZ_+$ is a double array of numbers satisfying the
following conditions:

\medskip

\begin{enumerate}
\item[$(i)$] There exist numbers $A_m$ such that for every $m\in\BZ_+$ and every $\delta>0$
\[
\sum_n a^\delta_{n,m}\leq A_m.
\]
\item[$(ii)$] For every $m\in\BZ_+$ 
\[
\lim_{\delta\to 0,\; N\to \infty} \sum_{n\geq N} a^\delta_{n,m}=0.
\]
\item[$(iii)$] $\sum_m A_m < \infty$.
\end{enumerate}

We claim that
\[
\lim_{\delta\to 0,\; N\to \infty} \sum_{n\geq N} \sum_m a^\delta_{n,m}=0.
\]
\end{lm}

\medskip

Choose $U_0=U_0(x_0)$ and $\Lambda=\lambda^{-1}$ according to Lemma 2.13 of \cite{KSSZ90}.
Then, by Poincar\'e recurrence, almost every point of $U_0$ infinitely many often returns to
$U_0$ and, on $U_0$ we have the uniform lower bound $\Lambda$ for the expansion rate of the
first return map restricted to {\it convex hypersurfaces}. Theorem \ref{strong_tail_bound} will be proved
once we check the hypotheses of Lemma \ref{arithmetic} for the numbers
\[
a^\delta_{n,m}=\delta^{-1}\cdot \mu_1\left(\tilde{U}^b_{n,m}(\delta)\right).
\]
In fact, by Lemma \ref{arithmetic} 
\[
\mu_1\left(\tilde{U}^b_\omega \right)\leq \sum_{n>F(\delta)} \sum_m 
\mu_1\left(\tilde{U}^b_{n,m}(\delta)\right)=o(\delta).
\]
To check the hypotheses of Lemma \ref{arithmetic} the following lemma will play an important role.

\medskip

\begin{lm} \label{disjoint}
For every fixed $m\in\BZ_+$, $0<n_1<n_2$ we have
\[
T^{n_1}\left(\tilde{U}^b_{n_1,m}\right)\inter T^{n_2}\left(\tilde{U}^b_{n_2,m}\right)
=\emptyset.
\]
\end{lm}

\begin{sprf}
Assume, on the contrary, that there exists a point
\[
w\in T^{n_1}\left(\tilde{U}^b_{n_1,m}\right)\inter T^{n_2}\left(\tilde{U}^b_{n_2,m}\right).
\]
Then, from the definition of the sets $\tilde{U}^b_{n,m}$, we have
\begin{equation} \label{inverses}
T^{-n_1}w,\; T^{-n_2}w\in U_0,
\end{equation}
and
\begin{equation} \label{kappas}
\kappa_{n_1,c_3\delta}(T^{-n_1}w),\; \kappa_{n_2,c_3\delta}(T^{-n_2}w)\in [\Lambda^m,\, \Lambda^{m+1}).
\end{equation}
But any $\tilde{\Sigma}$ permitted in the definition of $\kappa_{n_2,c_3\delta}(T^{-n_2}w)$
is also permitted in the definition of $\kappa_{n_1,c_3\delta}(T^{-n_1}w)$ and, moreover, due to 
(\ref{inverses}), the minimum expansion rate of $T^{n_2-n_1}$ on the manifold
$T^{n_1}\left(\tilde{\Sigma}\right)$ is at least $\Lambda$. This contradicts to
(\ref{kappas}), thus the lemma is proved.
\end{sprf} \qed

\medskip

By Lemma 5.3 and Lemma 4.10 of \cite{KSSZ90} we have
\[
\begin{aligned}
& \delta^{-1}\cdot \sum_n \mu_1\left(\tilde{U}^b_{n,m}\right)
=\delta^{-1}\cdot \sum_n \mu_1\left(T^n\tilde{U}^b_{n,m}\right)
=\delta^{-1}\cdot \mu_1\left(\bigcup_n T^n\tilde{U}^b_{n,m}\right) \\
& \leq\delta^{-1}\cdot \mu_1\left\{y\big|\; \exists n\;\; z_{\text{tub}}(y)<
\left(\kappa_{n,c_3\delta}(T^{-n}y)\right)^{-1}\cdot c_3\delta,\;
\kappa_{n,c_3\delta}(T^{-n}y)\in [\Lambda^m,\, \Lambda^{m+1})\right\} \\
& \leq\delta^{-1}\cdot \mu_1\left\{y\big|\; z_{\text{tub}}(y)<
\Lambda^{-m}\cdot c_3\delta\right\}\leq c_2c_3\Lambda^{-m}:=A_m,
\end{aligned}
\]
thus immediately implying that the hypotheses (i) and (iii) of Lemma \ref{arithmetic}
are satisfied with the numbers 
$a^\delta_{n,m}=\delta^{-1}\cdot \mu_1\left(\tilde{U}^b_{n,m}\right)$. The last remaining hypothesis (ii)
will follow from our

\begin{lm} \label{ordo}
For every fixed $m\in\BZ_+$ and every function $F(\delta)$ with $F(\delta)\nearrow \infty$ as
$\delta\searrow 0$, we have
\[
\lim_{\delta\to 0} \delta^{-1}\cdot
\sum_{n>F(\delta)} \mu_1\left(\tilde{U}^b_{n,m}(\delta)\right)=0.
\]
\end{lm}

\medskip

\begin{sprf}
We use the following, simple consequence of sufficiency: For every phase point $y\in\mathcal{S}_0$,
with a non-singular and sufficient backward orbit $\Phi^{(-\infty,0)}y$, it is true that
\begin{equation} \label{grows}
\lim_{n\to \infty} \kappa_{n,0}(T^{-n} y)=\infty.
\end{equation}
(This property is easily seen to be equivalent to the sufficiency of $\Phi^{(-\infty,0)}y$.)

Denote by $\mathcal{R}^*$ the set of all phase points $y\in\mathcal{S}_0$ for which the backward
orbit $\Phi^{(-\infty,0)}y$ is nonsingular and $\Phi^{(-\infty,-\tau(-y))}y$ is sufficient, i. e. the backward
orbit of $T^{-1} y\in\partial\mathbf{M}$ is sufficient. We define a smooth (actually, analytic) line bundle
$\mathcal{L}(\, .\, )$ over the manifold
\[
J=\left\{y\in\text{Int}(\mathbf{M})\big|\; Ty\in \mathcal{S}_0 \right\}
\]
as follows: For $y\in J$ let $n=n(y)=(z,\, w)=(z(y),\, w(y))$ be a smooth field of unit normal vectors
of $J$, and for $y=(q,\, v)\in J$ we take
\begin{equation} \label{bundle}
\begin{aligned}
\mathcal{L}(y)&=\left\{(q+tz(y),\, v)\big|\; t\geq 0 \text{ and } \forall \tau,\; 0\leq\tau\leq t,\;
(q+\tau z(y),\, v)\in\text{Int}(\mathbf{M})\right\} \\
& \cup \left\{(q+tz(y),\, v)\big|\; t\leq 0 \text{ and } \forall \tau,\; t\leq\tau\leq 0,\;
(q+\tau z(y),\, v)\in\text{Int}(\mathbf{M})\right\}.
\end{aligned}
\end{equation}
For $t_1\leq 0\leq t_2$ we introduce the notation
\[
\mathcal{L}^{t_2}_{t_1}(y)=\left\{(q+\tau z,\, v)\big|\; t_1\leq \tau \leq t_2\right\},
\]
provided that $\mathcal{L}^{t_2}_{t_1}(y)\subset \mathcal{L}(y)$. It is clear that the line bundle
$\mathcal{L}(y)$ ($y\in J$) is flow-invariant, i. e. 
\[
\Phi^s\left(\mathcal{L}^{t_2}_{t_1}(y)\right)=\mathcal{L}^{t_2}_{t_1}\left(\Phi^s(y)\right)
\]
for all $s$ with
\[
0\leq s<\min\left\{\tau(x)\big|\; x\in\mathcal{L}^{t_2}_{t_1}(y)\right\},
\]
and similar property holds true for the backward images of $\mathcal{L}^{t_2}_{t_1}(y)$.

In the sequel we will need the backward projection (i. e. the projection under $T^{-1}$)
line bundle $\mathcal{L}^*(y)$ ($y\in J\inter \partial\mathbf{M}=\mathcal{S}_{-1}$) of the line
bundle $\mathcal{L}(y)$ ($y\in J$), which is defined as follows:
\begin{equation} \label{l_star}
\begin{aligned}
\mathcal{L}^*(y)&=\big\{y'\in\partial\mathbf{M}\big|\; \exists t_1\, t_2,\, 0<t_1<\tau(y), \\
& 0<t_2<\tau(y'), \Phi^{t_2}(y')\in \mathcal{L}\left(\Phi^{t_1}(y)\right)\big\},
\end{aligned}
\end{equation}
for $y\in J\inter \partial\mathbf{M}=\mathcal{S}_{-1}$.

The definition of the sets $\tilde{U}^b_n(\delta)$ in (\ref{bad_tilde}) shows that for every phase point
$x\in \tilde{U}^b_n(\delta)$ there exists a number $\vep_1$, $0<\vep_1<\tau(T^n x)$, and there exists a
point $y\in J$ such that
\begin{equation} \label{sync}
\Phi^{\vep_1}(T^n x)\in\mathcal{L}(y),\; \text{dist}\left(\Phi^{\vep_1}(T^n x),\, y\right)
< \left(\kappa_{n,c_3\delta}(x)\right)^{-1}\cdot c_3\delta.
\end{equation}
Furthermore, for the phase point $y$, associated with $x$ in (\ref{sync}), it is true that
\begin{equation} \label{singular}
T y\in\mathcal{R}^*.
\end{equation}
This property is the motivation behind the definition of the set 
$\tilde{\mathcal{R}}=T^{-1}(\mathcal{R}^*)\inter \partial\mathbf{M}$. 

Now (\ref{sync}) and (\ref{singular}) say that for every point $x\in\tilde{U}^b_n(\delta)$
there exists a (unique) point $y\in\tilde{\mathcal{R}}$ such that
\begin{equation} \label{sync2}
T^n x\in\mathcal{L}^*(y),\; \text{dist}(T^n x,\, y)
< \left(\kappa_{n,c_3\delta}(x)\right)^{-1}\cdot c_3\delta.
\end{equation}

Let the number $m\in\BZ_+$ now be fixed. For every $y\in\tilde{\mathcal{R}}$ there exists an integer
$n_0=n_0(y,m)\in\BZ_+$ such that
\begin{equation} \label{kappa_big}
\kappa_{n_0,0}\left(T^{-n_0} y\right) > \Lambda^{m+2}.
\end{equation}
Then, thanks to the continuity of $\kappa_{n_0,0}\left(T^{-n_0} y\right)$ in $y$, for every
$y\in\tilde{\mathcal{R}}$ there exists an open neighborhood $V_1=V_1(y,m)\subset\partial\mathbf{M}$
of $y$ in $\partial\mathbf{M}$ such that for every $w\in V_1$ 
\begin{equation} \label{kappa_big2}
\kappa_{n_0,0}\left(T^{-n_0} w\right) > \Lambda^{m+2}.
\end{equation}
Also by continuity and by the definition $\kappa_{n,\delta}\left( y\right)$, for every 
$y\in\tilde{\mathcal{R}}$ there exists an open neighborhood $V_2=V_2(y,m)\subset\partial\mathbf{M}$
of $y$ in $\partial\mathbf{M}$ and a number $\delta_1=\delta_1(y,m)>0$ such that for every $\delta$,
$0<\delta<\delta_1$, and for every $w\in V_2$
\begin{equation} \label{kappa_big3}
\kappa_{n_0,\delta}\left(T^{-n_0} w\right) > \Lambda^{m+2}.
\end{equation}

Finally, since $\kappa_{n,\delta}(T^{-n} y)$ is non-decreasing in $n$, we also have that for every $\delta$,
$0<\delta<\delta_1$, for every $w\in V_2$, and for every $n\geq n_0$
\begin{equation} \label{kappa_big4}
\kappa_{n,\delta}\left(T^{-n} w\right) > \Lambda^{m+2}.
\end{equation}
Next we claim that for every $y\in\tilde{\mathcal{R}}$ and for every $\delta$, $0<\delta<c_3^{-1}\delta_1$,
\begin{equation} \label{image_disjoint}
\left[\bigcup_{n\geq n_0} T^n\left(\tilde{U}^b_{n,m}(\delta)\right)\right]
\inter V_2 = \emptyset.
\end{equation}
Indeed, 
\[
T^n\left(\tilde{U}^b_{n,m}(\delta)\right)\subset \big\{w\in\partial\mathbf{M}\big|\;
\kappa_{n,c_3\delta}(T^{-n} w)\in[\Lambda^m,\, \Lambda^{m+1}) \big\},
\]
thus, according to (\ref{kappa_big4}), the claim is true whenever $n\geq n_0$ and 
$c_3\delta<\delta_1$.

Due to the regularity of the hypersurface measure $\nu_0$ of 
$\mathcal{S}_{-1}=J\inter \partial\mathbf{M}$, for every $\eta>0$ there exists a compact set
$K_\eta\subset\tilde{\mathcal{R}}$ such that
\begin{equation} \label{K_big}
\nu_0\left(\tilde{\mathcal{R}}\setminus K_\eta\right) < \eta.
\end{equation}
Then, due to the compactness of $K_\eta$, one can choose a finite subset
\[
\{y_1,y_2,\dots,y_l\}\subset K_\eta
\]
such that
\begin{equation} \label{cover}
\bigcup_{i=1}^l V_2(y_i,\, m)\supset K_\eta,
\end{equation}
implying the existence of a $\delta_\eta > 0$ such that
\begin{equation} \label{cover2}
\bigcup_{i=1}^l V_2(y_i,\, m)\supset K_\eta^{[\delta_\eta]},
\end{equation}
where $A^{[\delta]}$ denotes the open $\delta$-neighborhood of a subset $A\subset\partial\mathbf{M}$
in $\partial\mathbf{M}$. We may assume that
\[
\delta_\eta < \min_{1\leq i\leq l} \delta_1(y_i,\, m).
\]
We take 
\[
N_\eta = \max_{1\leq i\leq l} n_0(y_i,\, m).
\]
Then, for every $\delta < c_3^{-1}\delta_\eta$ one has
\[
\left[\bigcup_{n\geq N_\eta} T^n\left(\tilde{U}^b_{n,m}(\delta)\right)\right]\inter
\left[\bigcup_{i=1}^l V_2(y_i,\, m)\right] = \emptyset,
\]
hence
\begin{equation} \label{disjoint2}
\left[\bigcup_{n\geq N_\eta} T^n\left(\tilde{U}^b_{n,m}(\delta)\right)\right]\inter
K_\eta^{[\delta_\eta]} = \emptyset.
\end{equation}

On the other hand, according to (\ref{sync2}), for $\delta < c_3^{-1}\delta_\eta$ and
$n\geq N_\eta$ one has
\begin{equation} \label{inclusion2}
T^n\left(\tilde{U}^b_{n,m}(\delta)\right) \subset \bigcup_{y\in\tilde{\mathcal{R}}} 
\left(\mathcal{L}^*(y)\right)^{[c_3\delta]},
\end{equation}
where 
\[
\left(\mathcal{L}^*(y)\right)^{[\vep]} = \big\{y'\in \mathcal{L}^*(y)\big|\; \text{dist}(y,\, y') 
< \vep \big\}.
\]
Consequently, whenever $c_3\delta < \delta_\eta$, we have
\begin{equation} \label{inclusion3}
\bigcup_{n\geq N_\eta}
T^n\left(\tilde{U}^b_{n,m}(\delta)\right) \subset \bigcup_{y\in\tilde{\mathcal{R}}\setminus K_\eta}
\left(\mathcal{L}^*(y)\right)^{[c_3\delta]}.
\end{equation}
By Lemma \ref{disjoint}
\begin{equation} \label{additive}
\sum_{n\geq N_\eta} \mu_1\left[T^n\left(\tilde{U}^b_{n,m}(\delta)\right)\right]
= \mu_1\left[\bigcup_{n\geq N_\eta} T^n\left(\tilde{U}^b_{n,m}(\delta)\right)\right]
\end{equation}
and, moreover, by (\ref{inclusion3}), by the smoothness of the bundle $\mathcal{L}^*(y)$,
and by (\ref{K_big}) the right-hand-side of \ref{additive} is at most $3c_3\eta\delta$. (Actually, here
the constant $3$ can be made arbitrarily close to $2$, but bigger than $2$, provided that
$\delta>0$ is small enough.) Here the multiplier $3c_3\eta$ of $\delta$ can be made arbitrarily small
by choosing the number $\eta>0$ small enough. Thus, we get that for any, arbitrarily small, given
constant $c_4>0$
\[
\sum_{n\geq N} \mu_1\left[\tilde{U}^b_{n,m}(\delta)\right] \leq c_4\delta,
\]
provided that $N$ is big enough and $\delta>0$ is small enough. This completes the proof of 
Lemma \ref{ordo} which, in turn, finishes the proof of Theorem \ref{strong_tail_bound} 
without using the Ansatz, thus our theorems \ref{main_result1} and \ref{main_result2}
are now also proved. \qed
\end{sprf}
\end{sprf}

\bigskip

\noindent
\textbf{Acknowledgements.} The author is deeply indebted to Nikolai I. Chernov and Imre P. T\'oth 
for their particularly useful remarks and suggestions.

\end{document}